# Analysis of Strongly Non-linear Oscillators by He's Improved Amplitude-Frequency Formulation

**R. Azami, D. D. Ganji[1], A. G. Davodi, H. Babazadeh**

*Department of Mechanical, Electrical and Civil Engineering, Babol University of Technology, Babol, Iran*

**Abstract**

In this work, we implement a relatively new analytical technique, the Improved Amplitude-Frequency Formulation (IAFF) method, approach for solving accurate approximate analytical solutions for strong nonlinear oscillators, which may contain high nonlinear terms.

This method can be used to obtain analytic and approximate solutions of different types of fractional differential equations applied in engineering mathematics. The solutions obtained are compared with those obtained by the Harmonic balance method (HBM) and Exact method, showing excellent agreement. We find that these attained solutions are not only with high degree of accuracy, but also uniformly valid in the whole solution domain which are so simple-to-do and effective.

***Keywords***: *Nonlinear oscillation, Improved Amplitude-frequency Formulation, period, exact solution, Harmonic balance method, fractional differential equation*

## 1. Introduction

We all know nonlinear functions are crucial points in engineering problems, so solving these equations are in the circle of most scientists and engineers' priority and requirements. Traditional perturbation methods which were methods for solving nonlinear equations have many shortcomings and they are not valid for strongly nonlinear equations. To overcome the shortcomings, many new techniques have appeared in open literature, such as, Delta-perturbation method [1, 2], Variational Iteration method (VIM) [3-7], Homotopy Perturbation method (HPM) [10-18], Bookkeeping Parameter Perturbation method [19], He's Energy Balance method (EBM) [20-21] and the Improved Amplitude-frequency Formulation (IAFF) [22] and other methods for different problems [29- 36]. In IAFF method the angular frequency can be readily obtained. The results are valid not only for weakly nonlinear systems, but also for strongly nonlinear ones.

## 2. Basic Idea of Improved Amplitude-frequency Formulation

We consider a generalized nonlinear oscillator in the form:

$$u'' + f(u) = 0, \quad u(0) = A, \quad u'(0) = 0. \tag{1}$$

We use two following trial functions:

---

[1] Corresponding author:
E-mail: ddg_davood@yahoo.com, Tel.: +98 111 3234205 (D.D. Ganji)

$$u_1(t) = A \cos \omega_1 t \tag{2}$$

and

$$u_2(t) = A \cos \omega_2 t \tag{3}$$

The residuals are

$$R_1(t) = -A \cos \omega_1 t + f(\cos \omega_1 t) \tag{4}$$

and

$$R_2(\omega t) = -A \omega_2^2 \cos \omega_2 t + f(\cos \omega_2 t) \tag{5}$$

The original Frequency-amplitude formulation reads [23-26]:

$$\omega^2 = \frac{\omega_1^2 R_2 - \omega_2^2 R_1}{R_2 - R_1} \tag{6}$$

He used the following formulation [23-26] and Geng and Cai improved the formulation by choosing another location point [27].

$$\omega^2 = \frac{\omega_1^2 R_2(\omega_2 t = 0) - \omega_2^2 R_1(\omega_1 t = 0)}{R_2 - R_1} \tag{7}$$

This is the improved form by Geng and Cai.

$$\omega^2 = \frac{\omega_1^2 R_2(\omega_2 t = \pi/3) - \omega_2^2 R_1(\omega_1 t = \pi/3)}{R_2 - R_1} \tag{8}$$

The point is: $\cos \omega_1 t = \cos \omega_2 t = k$

Substituting the obtained ω in to $u(t) = A\cos \omega t$, we can obtain the constant $k$ in $\omega^2$ equation in order to have the frequency without irrelevant parameter.

To improve its accuracy, we can use the following trial function when they are required.

$$u_1(t) = \sum_{i=1}^{m} A_i \cos \omega_i t \quad and \quad u_2(t) = \sum_{i=1}^{m} A_i \cos ?_i t \tag{9}$$

or

$$u_1(t) = \frac{\sum_{i=1}^{m} A_i \cos \omega_i t}{\sum_{j=1}^{n} B_j \cos \omega_j t} \quad and \quad u_2(t) = \frac{\sum_{i=1}^{m} A_i \cos ?_{ii} t}{\sum_{j=1}^{n} B_j \cos ?_{ij} t} \tag{10}$$

But in most cases because of the sufficient accuracy trial functions are as follow and just the first term:

$$u_1(t) = A\cos(t) \quad and \quad u_2(t) = A\cos(\omega t) + (A-a)\cos(3\omega t) \tag{11}$$

And

$$u_1(t) = A \cos t \quad and \quad u_2(t) = \frac{A(1+c) \cos \omega t}{1 + c \cos 2\omega t} \tag{12}$$

Where $a$ and $c$ are unknown constants. In addition we can set $\cos t = k$ in $u_1$ and $\cos \omega t = k$ in $u_2$.

### 3. Application

To show the efficiency and accuracy of the method, in the following we solve some examples by IAFF and compare the results with those of HBM and Exact methods.

### 3.1. Example 1:

Consider the motion of a particle of mass *m* attached to the center of a stretched elastic wire [37] of coefficient of stiffness equal to *k*. The length of elastic wire when any force is applied to it is *2a*. We assume that the movement of particle is one-dimensional and this is constrained to move only in the horizontal *x* direction.

As we can see in Fig. 1, ends of wire are fixed a distance *2d* a part. Length *d* can be major or equal to *a*. If $d = a$, the wire is not stretched for $x = 0$, and there is no tension in each part of it. However, if $d > a$ the wire is stretched for $x = 0$, and the tension in each part of the wire is $k(d–a)$. The equation of motion is given by the following nonlinear differential equation [38]:

$$m\frac{d^2x}{dt^2} + 2kx - \frac{2kax}{\sqrt{d^2+x^2}} = 0 \tag{13}$$

With initial conditions:

$$x(0) = A \ , \ \frac{dx}{dt}(0) = 0 \tag{14}$$

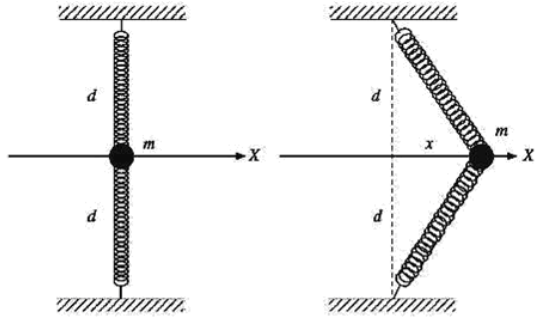

Figure 1: Mass attached to a stretched wire.

*u* and *t* can be constructed as follows:

$$u = dx \ , t = \sqrt{\frac{2k}{m}} \tag{15}$$

And:

$$\frac{d^2u}{dt^2} + u - \frac{\lambda u}{\sqrt{1+u^2}} = 0 \ , 0 < \lambda \leq 1 \tag{16}$$

With initial conditions:

$$u(0) = A \ , \ \frac{du}{dt}(0) = 0 \tag{17}$$

For small values of *A* we can write:

$$\frac{1}{\sqrt{1+u^2}} = \left(1 - \frac{1}{2}u^2\right) \tag{18}$$

We can write Eq. (16) in the form of:

$$u'' + (1-\lambda)u + \frac{1}{2}\lambda u^3 = 0, \quad u(0) = A, u'(0) = 0 \tag{19}$$

With initial conditions of:

$$u(0) = A, u'(0) = 0. \tag{20}$$

The main equation and its conditions are:

$$u_1(t) = A \cos t \tag{21}$$

and

$$u_2(t) = A \cos 2t \tag{22}$$

Considering $u_1(t) = u_1(t) = A \cos t = A \cos 2t = k$, the residual equations are, respectively:

$$R_1(t) = -\lambda A k + \frac{1}{2}\lambda A^3 k^3 \tag{23}$$

and

$$R_2(2t) = -3Ak - \lambda A k + \frac{1}{2}\lambda A^3 k^3 \tag{24}$$

Considering to $\cos \omega_1 t = \cos \omega_2 t = k$, we have:

$$\omega^2 = \frac{\omega_1^2 R_2 - \omega_2^2 R_1}{R_2 - R_1} = 1 - \lambda + \frac{1}{2}\lambda A^2 k^2 \tag{25}$$

$$\omega = \sqrt{1 - \lambda + \frac{1}{2}\lambda A^2 k^2} \tag{26}$$

We can rewrite $u(t) = A \cos \omega t$, in the form:

$$u(t) = A\cos\left[\sqrt{1 - \lambda + \frac{1}{2}\lambda A^2 k^2}\, t\right] \tag{27}$$

In view of the approximate solution, we can rewrite the main equation in the form:

$$u'' + \left((1-\lambda) + \frac{1}{2}\lambda A^2 k^2\right)u = \frac{1}{2}\lambda A^2 k^2 u - \frac{1}{2}\lambda u^3 \tag{28}$$

If by any chance $u(t) = A\cos\left[\sqrt{1 - \lambda + \frac{1}{2}\lambda A^2 k^2}\, t\right]$ is the exact solution, then the right side of Eq. (28) is vanishing completely. Considering to our approach which is just an approximation one, we set:

$$B = \int_0^{T/4} \left(\frac{1}{2}\lambda A^2 k^2 u - \frac{1}{2}\lambda u^3\right) \cos \omega t\, dt = 0, \quad T = \frac{2\pi}{\omega}. \tag{29}$$

$$B_1 = \int_0^{T/4} \left(\frac{1}{2}\lambda A^2 k^2 A \cos \omega t\right) \cos \omega t\, dt = \frac{1}{2}\lambda A^3 k^2 \frac{\pi}{4} \tag{30}$$

$$B_2 = \int_0^{T/4} \left(-\frac{1}{2}\lambda(A\cos\omega t)^3\right)\cos\omega t\, dt = -\frac{1}{2}\lambda A^3 \frac{3\pi}{16} \tag{31}$$

$$B = B_1 + B_2 \rightarrow \frac{1}{2}\lambda A^3 k^2 \frac{\pi}{4} - \frac{1}{2}\lambda A^3 \frac{3\pi}{16} = 0 \tag{32}$$

Solving and simplifying Eq. (32), we have:

$$k^2 = \frac{3}{4} \tag{33}$$

So

$$\omega = \sqrt{1 - \lambda + \frac{3}{8}\lambda A^2} \tag{34}$$

In order to compare with other methods, we brought Beléndez obtained Harmonic balance solutions for this example [38]:

$$\omega_{HBM} = \sqrt{1 - \frac{\lambda}{\sqrt{f(A)}}}, \sqrt{f(A)} = 1 + \frac{3}{4}A^2 + \frac{3}{64}A^4 + \frac{16}{512}A^6 + \ldots \tag{35}$$

Figure 2 shows Comparison of the Harmonic balance solution with the Improved Amplitude-frequency Formulation when we substitute $A = 0.3$.

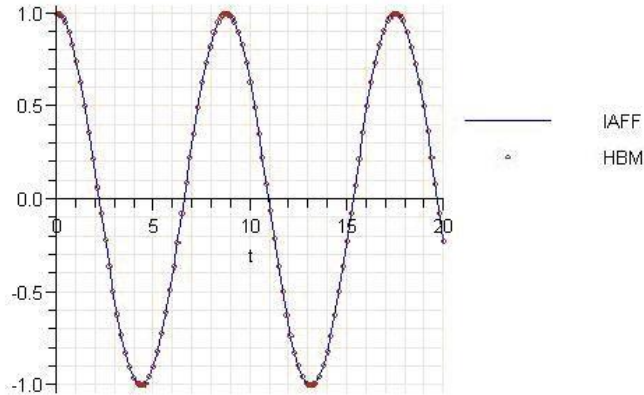

Figure 2: comparison of IAFF and HBM for (A=0.3).

Table1: For a constant $\lambda = 0.5$ the table shows the comparison of IAFF result with HBM results.

| A | ω$_{HBM}$ | ω$_{IAFF}$ |
|---|---|---|
| 0.02 | 0.7072 | 0.70715 |
| 0.04 | 0.7073 | 0.70731 |
| 0.1 | 0.7084 | 0.70843 |
| 0.16 | 0.7104 | 0.71049 |
| 0.2 | 0.7123 | 0.71239 |
| 0.3 | 0.71833 | 0.71894 |
| 0.4 | 0.72616 | 0.72801 |
| 0.5 | 0.73531 | 0.73950 |
| 0.6 | 0.74531 | 0.75332 |
| 0.7 | 0.75573 | 0.76933 |
| 0.8 | 0.76632 | 0.7874 |

### 3.2. Example 2

This example corresponds to the following oscillator:

$$u'' + u^{\frac{3}{4}} = 0 \tag{36}$$

With initial conditions of:

$$u(0) = A, u'(0) = 0. \tag{37}$$

The main equation and its conditions are:

$$u_1(t) = A \cos t \tag{38}$$

and

$$u_2(t) = A \cos 2t \tag{39}$$

Considering $u_1(t) = u_1(t) = A \cos t = A \cos 2t = k$, the residual equations are, respectively:

$$R_1(t) = -Ak + A^{\frac{3}{4}} k^{\frac{3}{4}} \tag{40}$$

and

$$R_2(2t) = -4Ak + A^{\frac{3}{4}} k^{\frac{3}{4}} \tag{41}$$

Considering to $\cos \omega_1 t = \cos \omega_2 t = k$, we have:

$$\omega^2 = \frac{\omega_1^2 R_2 - \omega_2^2 R_1}{R_2 - R_1} = \frac{1}{A^{\frac{1}{4}} k^{\frac{1}{4}}} \tag{42}$$

$$\omega = \frac{1}{A^{\frac{1}{2}} k^{\frac{1}{2}}} \tag{43}$$

We can rewrite $u(t) = A \cos \omega t$, in the form:

$$u(t) = A\cos\left[\left(A^{\frac{-1}{2}} k^{\frac{-1}{2}}\right)t\right] \tag{44}$$

In view of the approximate solution, we can rewrite the main equation in the form:

$$u'' + \frac{1}{A^{\frac{1}{4}} k^{\frac{1}{4}}} u = \frac{1}{A^{\frac{1}{4}} k^{\frac{1}{4}}} u - u^{\frac{3}{4}} \tag{45}$$

If by any chance Eq. (44) is the exact solution, then the right side of Eq. (45) is vanishing completely. Considering to our approach which is just an approximation one, we set:

$$B = \int_0^{T/4} \left( \frac{1}{A^{\frac{1}{4}} k^{\frac{1}{4}}} u - u^{\frac{3}{4}} \right) \cos \omega t \, dt = 0, \quad T = \frac{2\pi}{\omega}. \tag{46}$$

$$B_1 = \int_0^{T/4} \left( \frac{1}{A^{\frac{1}{4}} k^{\frac{1}{4}}} A \cos \omega t \right) \cos \omega t \, dt = A^{\frac{3}{4}} k^{-\frac{1}{4}} \frac{\pi}{4} \tag{47}$$

$$B_2 = \int_0^{T/4}\left(-(A\cos\omega t)^{\frac{3}{4}}\right)\cos\omega t\, dt = -A^{\frac{3}{4}}\frac{\frac{3}{14}\pi^{\frac{3}{2}}\csc\left(\frac{3\pi}{8}\right)}{\Gamma\left(\frac{5}{8}\right)\Gamma(\frac{7}{8})} \tag{48}$$

$$B = B_1 + B_2 \;\rightarrow\; A^{\frac{3}{4}}k^{-\frac{1}{4}}\frac{\pi}{4} - A^{\frac{3}{4}}\frac{\frac{3}{14}\pi^{\frac{3}{2}}\csc\left(\frac{3\pi}{8}\right)}{\Gamma\left(\frac{5}{8}\right)\Gamma(\frac{7}{8})} = 0 \tag{49}$$

Solving and simplifying Eq. (41), we have:

$$k = \left(\frac{\Gamma\left(\frac{5}{8}\right)\Gamma\left(\frac{7}{8}\right)}{\frac{12}{14}\pi^{\frac{1}{2}}\csc\left(\frac{3\pi}{8}\right)}\right) = 0.9505641275 \;\rightarrow\; k\cong 0.95 \;\rightarrow\; k^{\frac{1}{2}}=0.97 \tag{50}$$

So

$$\omega = \frac{1}{A^{\frac{1}{2}}k^{\frac{1}{2}}} = \frac{1}{0.97 A^{\frac{1}{2}}} = \frac{1.030927}{A^{\frac{1}{2}}} \tag{51}$$

For $A=1$, we obtain $\omega_{IAFF} = 1.030927$, in comparison with exact solution $\omega_{ex} = 1.024957$ [39], and the error is just 0.57%. In this example we can obtain the frequency for different values of amplitude.

Figure 3 shows comparison of IAFF frequency with exact frequency when $A=1$.

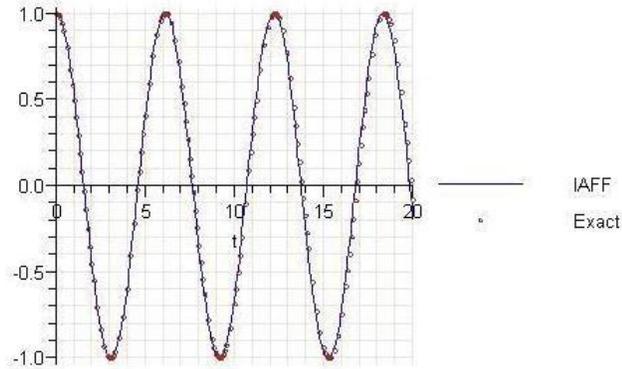

Figure 3: Comparison of IAFF frequency with exact frequency ($A=1$).

### 3.3. Example 3

Consider the following nonlinear oscillator with discontinuities [40,41]:

$$u'' + u + \varepsilon\, u^2 sgn(u) + u^3 = 0 \tag{52}$$

Where $sgn(u)$ is the sign function, the result equals to *+1* if *u>0*, *0* if *u = 0* and *-1* if *u<0*.

With initial conditions of:

$$u(0) = A, u'(0) = 0 \tag{53}$$

The main equations are:
$$u_1(t) = A \cos t \tag{54}$$

And
$$u_2(t) = A \cos 2t \tag{55}$$

If $u > 0$:
$$u'' + u + \varepsilon u^2 + u^3 = 0 \tag{56}$$

Respectively the residual equations are:
$$R_1(t) = -A \cos t + A \cos t + \varepsilon A^2 \cos^2 t + A^3 \cos^3 t \tag{57}$$

And
$$R_2(2t) = -4A \cos 2t + A \cos 2t + \varepsilon A^2 \cos^2 2t + A^3 \cos^3 2t \tag{58}$$

Considering to $\cos \omega_1 t = \cos \omega_2 t = k$, we have:

$$\omega^2 = \frac{\omega_1^2 R_2 - \omega_2^2 R_1}{R_2 - R_1} = 1 + \varepsilon A k + A^2 k^2 \rightarrow \omega = \sqrt{1 + \varepsilon A k + A^2 k^2} \tag{59}$$

We can rewrite $u(t) = A \cos \omega t$ in the form:

$$u(t) = A \cos[(1 + \varepsilon A k + A^2 k^2)^{\frac{1}{2}} t] \tag{60}$$

In view of the approximate solution, we can rewrite the main equation in the form:
$$u'' + (1 + \varepsilon A k + A^2 k^2) u = -\varepsilon u^2 - u^3 + \varepsilon A k u + A^2 k^2 u \tag{61}$$

If by any chance Eq. (60) is the exact solution, then the right side of Eq. (61) is vanishing completely. Considering to our approach which is just an approximation one, we set:

$$B = \int_0^{T/4} (-\varepsilon u^2 - u^3 + \varepsilon A k u + A^2 k^2 u) \cos \omega t \, dt = 0, \quad T = \frac{2\pi}{\omega}. \tag{62}$$

$$\rightarrow -\varepsilon \frac{A^2}{6} - \frac{3\pi}{16} A^3 + A \frac{\pi}{4} (\varepsilon A k + A^2 k^2) = 0 \tag{63}$$

Solving Eq. (63), we have:

$$k = \frac{-\pi \varepsilon \mp \sqrt{(\pi \varepsilon)^2 + 4A\pi(\frac{2}{3}\varepsilon + \frac{3}{4} A\pi)}}{2A\pi} \tag{64}$$

$$\omega = \sqrt{1 + \varepsilon A k + A^2 k^2} \tag{65}$$

If $u > 0$:
$$u'' + u - \varepsilon u^2 + u^3 = 0 \tag{66}$$

Respectively the residual equations are:
$$R_1(t) = -A \cos t + A \cos t - \varepsilon A^2 \cos^2 t + A^3 \cos^3 t \tag{67}$$

And
$$R_2(t) = -4A \cos 2t + A \cos 2t - \varepsilon A^2 \cos^2 2t + A^3 \cos^3 2t \tag{68}$$

Considering to $cos\ \omega_1 t = cos\ \omega_2 t = k$, we have:

$$\omega^2 = \frac{\omega_1^2 R_2 - \omega_2^2 R_1}{R_2 - R_1} = 1 - \varepsilon Ak + A^2 k^2 \quad \omega = \sqrt{1 - \varepsilon Ak + A^2 k^2} \tag{69}$$

We can rewrite $u(t) = A\ cos\ \omega t$ in the form:

$$u(t) = Acos\ [\,(1 - \varepsilon Ak + A^2 k^2)^{\frac{1}{2}} t\,] \tag{70}$$

In view of the approximate solution, we can rewrite the main equation in the form:

$$u'' + (1 - \varepsilon Ak + A^2 k^2)u = \varepsilon u^2 - u^3 - \varepsilon Aku + A^2 k^2 u \tag{71}$$

If by any chance Eq. (70) is the exact solution, then the right side of Eq. (71) is vanishing completely. Considering to our approach which is just an approximation one, we set:

$$B = \int_0^{T/4} (\varepsilon u^2 - u^3 - \varepsilon Aku + A^2 k^2 u)\ cos\ \omega t\ dt = 0, \quad T = \frac{2\pi}{\omega}. \tag{72}$$

$$\rightarrow \varepsilon \frac{A^2}{6} - \frac{3\pi}{16} A^3 + A \frac{\pi}{4} \left(-\varepsilon Ak + A^2 k^2\right) = 0 \tag{73}$$

Solving Eq. (73), we have:

$$k = \frac{\pi \varepsilon \mp \sqrt{(\pi \varepsilon)^2 + 4 A \pi (\frac{2}{3}\varepsilon - \frac{3}{4} A \pi)}}{2 A \pi} \tag{74}$$

$$\omega = \sqrt{1 - \varepsilon Ak + A^2 k^2} \tag{75}$$

We obtain the approximate period $T$:

$$T = \frac{T_1 + T_2}{2}$$

In order to compare with exact solution [40]:

$$T_e = \int_0^A \frac{2dx}{\sqrt{A^2 - x^2 + \frac{2}{3}\varepsilon(A^3 - x^3) + \frac{1}{2}(A^4 - x^4)}} + \int_0^A \frac{2dx}{\sqrt{A^2 - x^2 - \frac{2}{3}\varepsilon(A^3 - x^3) + \frac{1}{2}(A^4 - x^4)}} \tag{76}$$

Figure 4 shows Comparison of IAFF period with exact period when $\varepsilon = 1$ and $A = 5$.

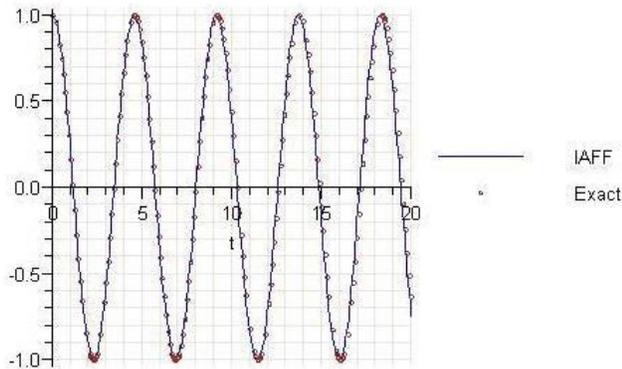

Figure 4: Comparison of IAFF period with exact period. $(\varepsilon = 1, A = 5)$

Table2: Comparison of IAFF period with exact period ($\varepsilon = 1$)

| A | IAFF period | Exact period | Error (%) |
|---|---|---|---|
| 0.1 | 6.2176870 | 6.28559 | 1.0 |
| 0.2 | 6.1558638 | 6.28794 | 2.1 |
| 0.4 | 6.0603034 | 6.20413 | 2.3 |
| 0.5 | 6.0399889 | 6.05759 | 0.29 |
| 5 | 1.3773142 | 1.36965 | 0.55 |
| 10 | 0.7108746 | 0.71463 | 0.52 |
| 100 | 0.0724447 | 0.07391 | 1.9 |

**4. Conclusions**

He's improved Amplitude-frequency Formulation, was applied to nonlinear oscillators which are useful in so many branches of sciences such as: fluid mechanics, electromagnetic and waves, telecommunication, civil and its structures and all so-called majors' applications and etc… The Amplitude-frequency Formulation is a well-established method for the analysis of nonlinear systems, can be easily extended to any nonlinear equation. We demonstrated the accuracy and efficiency of the presenting method with some strong nonlinear problems. Improved Amplitude-frequency Formulation provides an easy and direct procedure for determining approximations to the periodic solutions. Eventually, this paper suggests to readers to apply improved Amplitude-frequency Formulation method for solving nonlinear oscillations because of its accuracy, reliability and simplicity.